\documentclass[12pt]{article}
\usepackage{amsmath}
\usepackage{amssymb}
\usepackage{euscript}
\usepackage[latin1]{inputenc}

\normalbaselines
\newcommand{\qed}{\sqcap\kern-8pt\sqcup}

\newenvironment{proof}
      {\par\noindent{\it Proof\/: }\nopagebreak\normalsize}%
                                  {\linebreak[2]\hspace*{\fill}$\qed$\ifdim\lastskip<12pt
       \removelastskip \penalty-200  \vskip12pt  \fi}

\font\frten=eufm10 at 12pt
\font\freight=eufm10
\font\frsix=eufm8
\newfam\frfam\textfont\frfam=\frten
\scriptfont\frfam=\freight
\scriptscriptfont\frfam=\frsix

\newcommand{\CC}{\mathbb{C}}
\newcommand{\PP}{\mathbb{P}}
\newcommand{\RR}{\mathbb{R}}

\newcommand{\ZZ}{\mathbb{Z}}

\newtheorem{thm}{Theorem}[section]
\newtheorem{prop}[thm]{Proposition}
\newtheorem{cor}[thm]{Corollary}
\newtheorem{lem}[thm]{Lemma}
\newtheorem{defn}[thm]{Definition}

\newtheorem{rema}[thm]{Remark}

\def \O{{\cal O}}

\def \Div{{\rm Div}}
\def \Gal{{\rm Gal}}

\def \Aut{{\rm Aut}}
\def \Int{{\rm Int}}
\def \Stab{{\rm Stab}}

\begin{document}

\title{Fixed points of automorphisms of real algebraic curves.
\footnote{Work supported by the European Community's 
Human Potential Programme 
under contract HPRN-CT-2001-00271, RAAG.}}

\author{Jean-Philippe Monnier\\
       {\small D\'epartement de Mathématiques, Universit\'e d'Angers,}\\
{\small 2, Bd. Lavoisier, 49045 Angers cedex 01, France}\\
{\small e-mail: monnier@tonton.univ-angers.fr}}
\date{}
\maketitle

\begin{abstract}
{We bound the number of fixed points of an
automorphism of a real curve in terms of the genus and the number of
connected components of the real part of the curve. Using this
bound, we derive some consequences concerning the maximum order of
an automorphism and the maximum order of 
an abelian group of automorphisms of a real curve. We also bound the
full group of automorphisms of a real hyperelliptic curve.}
\end{abstract}

\section*{Introduction}

In this note, a real algebraic curve $X$ is a proper geometrically
integral scheme over $\RR$ of dimension $1$. Let $g$ denote the genus
of $X$; throughout the paper we assume $g\geq 2$.
A closed point $P$ of $X$ 
will be called a real point if the residue field at $P$ is $\RR$, and 
a non-real point if the residue field at $P$ is $\CC$.
The set of real points $X(\RR)$
will always be assumed to be non empty.
It decomposes into 
finitely many connected components, whose number will be denoted by $s$.
By Harnack's Theorem we know that $s\leq g+1$. 
If $X$ has $g+1-k$ real connected components, we will say that
$X$ is an $(M-k)$-curve. 
We will say that $X$ has many real components if $s\geq g$
(see \cite{Hu}).
Topologically, each semi-algebraic connected component of $X(\RR)$ 
is  homotopy equivalent
to a circle.

An automorphism $\varphi$ of $X$ is an isomorphism of schemes of $X$ with
itself. Seeing $X$ as a compact Klein surface, May proved in
\cite{Ma1} that the
order of a group of automorphisms of $X$ is bounded above by
$12(g-1)$.
Moreover, he also proved that the maximum possible order of an
automorphism of $X$ is $2g+2$ \cite{Ma2}.

In this paper, we will bound the number of fixed points of an
automorphism of a real curve in terms of $s$. Using this
bound, we will derive some consequences concerning the maximum order of
an automorphism and the maximum order of 
an abelian group of automorphisms of a real curve. We will also bound the
full group of automorphisms of a real hyperelliptic curve.

\section{ Preliminaries}

We recall here some classical concepts and more notation that we will be using 
throughout this paper.

Let $X$ be a real curve. We will denote by $X_{\CC}$ the ground field
extension of $X$ to $\CC$. 
The group $\Div(X)$ (resp. $\Div(X_{\CC})$) of divisors on $X$ (resp. $X_{\CC}$)
is the free abelian group generated by the 
closed points of $X$ (resp. $X_{\CC}$). The Galois group $\Gal(\CC /\RR)$
acts on the complex variety $X_{\CC}$ and also on $\Div(X_{\CC})$. 
We will always indicate this action by a bar.
If $P$ is a non-real point of $X$, identifying $\Div(X)$
and $\Div(X_{\CC})^{\Gal(\CC /\RR)}$, then
$P=Q+\bar{Q}$ 
with $Q$ a closed point of $X_{\CC}$. 

Let $\Aut(X)$ (resp. $\Aut(X_{\CC})$) denote the group of automorphisms
of $X$ (resp. $X_{\CC}$). If $\varphi\in \Aut(X)$ then $\varphi$
extends to an 
automorphism 
$\varphi_{\CC}$ of $X_{\CC}$ such that
$\overline{\varphi_{\CC} (Q)}=\varphi_{\CC }(\bar{Q})$ for
any closed point $Q$ of $X_{\CC}$. We will denote by $\mu (\varphi_{\CC})$
(resp. $\mu_{\RR} (\varphi)$)
the number of closed (resp. real closed) 
fixed points of $\varphi_{\CC }$ (resp. $\varphi$).

Let $G$ be a group and let $a$ be an element of $G$. We denote by
$|G|$ (resp. $|a|$) the order of $G$ (resp. $a$). 

Let $\varphi\in \Aut (X)$. The image of a connected component of $X
(\RR )$ through $\varphi$
is again a connected component of $X(\RR )$.
Let $\Sigma_s$ denote the symmetric group corresponding to the
group of permutations of the set of connected components of $X (\RR
)$. We will denote by $\sigma (\varphi)$ the permutation induced by $\varphi$.

Assume $G$ is a subgroup of $\Aut (X)$ of order $N$. 
We denote by $G_{\CC}$ the set $\{\varphi_{\CC
},\,\, \varphi\in G\}$. Then $G_{\CC }$ is a subgroup of $\Aut(X_{\CC})$ and 
$|G|=|G_{\CC }|$ since clearly $|\varphi|=|\varphi_{\CC }|$
if $\varphi\in \Aut(X)$. The quotient space $X/G$ is a real algebraic curve
of genus $g'$
and with above notation, $(X/G)_{\CC} =X_{\CC}/G_{\CC}$. We denote by
$\pi$ (resp. $\pi_{\CC}$) the morphism of real (resp. complex)
algebraic curves $X\rightarrow X/G$ (resp. $X_{\CC }\rightarrow
(X/G)_{\CC}$). We know that these maps are of degree $N$.
The map $\pi_{\CC}$ is ramified only at the fixed points of elements
of $G_{\CC}\setminus \{ Id\}$
and the ramification index $e_P$ at a closed point $P$ of $X_{\CC}$
verifies $e_P =|\Stab (P)|$, with $\Stab (P)$ the stabilizer
subgroup of $P$ in $G_{\CC}$. We say that two closed points $P, P'$ of
$X_{\CC }$ are equivalent under $G_{\CC}$ 
if there exists $\varphi$ in $G_{\CC}$ such that
$\varphi(P)=P'$. Using group theory, there are $\frac{N}{e_P}$ 
distinct points on
$X_{\CC}$ equivalent under $G_{\CC}$ to $P$. Consequently the number
$w(\pi_{\CC} )$ of branch points of $\pi_{\CC}$ corresponds to the
maximal number of inequivalent fixed points of elements of
$G_{\CC}\setminus \{ Id\}$. Let $Q$ denote a branch point of
$\pi_{\CC}$, we denote by $e(Q)$ the ramification index $e_P$ of any
ramification point $P$ over $Q$. Considering now the map $\pi
:X\rightarrow X/G$, we denote by $w_{\RR} (\pi )$ (resp. $w_{\CC} (\pi )$) the
number of real (resp. non-real) branch points of $\pi$.
Clearly $w_{\RR} (\pi )+2w_{\CC} (\pi )=w(\pi_{\CC} )$. If $Q$ is a real
branch point of $\pi $ then $Q$ corresponds to a branch point $Q'$ of
$\pi_{\CC}$ and $e(Q)=e(Q')$. If $Q$ is a non-real
branch point of $\pi $ then $Q$ corresponds to two conjugate branch
points 
$Q'$ and $\bar{Q'}$ of
$\pi_{\CC}$ and $e(Q)=e(Q')=e(\bar{Q'})$. Let $g'$ denote
the genus of $X/G$. The Riemann-Hurwitz relation now reads \cite[p. 243]{Fa-K}
\begin{equation}
\label{equ1}
2g-2=N(2g' -2) +N\sum_{i=1}^{w(\pi_{\CC})} (1-\frac{1}{e(Q_i)})
\end{equation}
where $Q_1 ,\ldots ,Q_{w(\pi_{\CC})}$ are the branch points of
$\pi_{\CC }$. Equivalently, we have
\begin{equation}
\label{equ2}
2g-2=N(2g' -2) +N\sum_{i=1}^{w_{\RR} (\pi)} (1-\frac{1}{e(Q_i)})
+2N\sum_{i=1}^{w_{\CC} (\pi)} (1-\frac{1}{e(Q'_i)})
\end{equation}
where $Q_1 ,\ldots ,Q_{w_{\RR} (\pi)}$ (resp. $Q'_1 ,\ldots ,Q'_{w_{\CC}
  (\pi)}$) are the real (resp. non-real) branch points of $\pi$.
If moreover $G=\langle\varphi\rangle$ is cyclic, then we have \cite[p. 245]{Fa-K}
\begin{equation}
\label{equ3}
2g-2=|\varphi|(2g' -2) +\sum_{i=1}^{|\varphi|-1} \mu (\varphi_{\CC}^i)
\end{equation}

\section{Real and non-real fixed points of  automorphisms of real curves}

 Before the study of the fixed points set of an automorphism of $X$ we
 need to state basic facts about automorphisms of real curves.

\begin{prop} 
\label{basic1}
Let $G$ be a subgroup of $\Aut(X)$. We denote by $X'$ the
  quotient real curve $X/G$.
\begin{description}
\item[({\cal{i}})] If $Q$ is a real point of $X'$ then either
  $\pi^{-1} (Q)$ is totally real or $\pi^{-1} (Q)$ is totally non-real.
\item[({\cal{ii}})] If $Q$ is a non-real point of $X'$ then 
  $\pi^{-1} (Q)$ is totally non-real.
\item[({\cal{iii}})] If $P$ is a real ramification point of $\pi$ then
  $e_P =2$.
\item[({\cal{iv}})] The image by $\pi$ of a connected component $C$ of
  $X(\RR)$ is a connected component of $X' (\RR)$ if and only if $C$
  does not contain any real ramification point of $\pi$. If not $\pi
  (C)$ is a compact connected semi-algebraic subset of a connected
  component $C'$ of $X'(\RR )$ corresponding topologically to a closed
  interval of $C'$, in this situation $\pi (C)$ contains exactly two real
  branch points of $\pi$ corresponding topologically to the end-points
  of $\pi (C)$. 
\item[({\cal{v}})] The number of real
  branch points of $\pi$ with real fibers is even on any connected component
  of $X'(\RR)$.
\item[({\cal{vi}})] The number of real
  ramification points of $\pi$ is even on any connected component
  of $X(\RR)$.
\end{description}
\end{prop}

\begin{proof} 
To prove statements ({\cal{i}}) and ({\cal{ii}}), it is sufficient to
make the observation that
an automorphism of $X$ maps a real (resp. non-real)
point onto a real (resp. non-real) point.

Let $P$ be a real ramification point of $\pi$. Then $Q=\pi (P)$ is a
real branch point of $\pi$. Let $t\in \O_Q$ be a local parameter. We
may consider $t$ as an element of $\O_P$ via $\pi^{\#}: \O_Q\rightarrow
\O_P$ and we recall that the ramification index $e_P$ is defined as
$v_P (t)$ where $v_P$ is the valuation associated to $\O_P$. We have
$e_P\geq 2$ since $P$ is a ramification point of $\pi$.
If $u$ is a
local parameter of $\O_P$ then $t=u^{e_P }a$ with $a\in
\O_P^*$. Consequently, $\pi$ is given locally at $P$ by $u\mapsto
u^{e_P}$. If $e_P >2$ then the fiber above a real point sufficiently
near $Q$ is not totally real or totally non-real contradicting statement
({\cal{i}}). Hence ({\cal{iii}}).

Statement ({\cal{iv}}) implies statements ({\cal{v}}) and ({\cal{vi}})
since $X'$ is a quotient space and, if $P$ is a real
ramification point of $\pi$ then any point in the fiber $\pi^{-1} (\pi
(P))$ is a real ramification point of $\pi$ with ramification index
equal to
$e_P =2$.

For ({\cal{iv}}), let $P$ be a real ramification point of $\pi$. Then
$Q=\pi (P)$ is a
real branch point of $\pi$. Let $C$ (resp. $C'$) denote the connected
component of $X(\RR )$ (resp. $X'(\RR)$) containing $P$ (resp. $Q$). 
Since $e_P =2$ is even, $C$ is clearly on one side of the fiber
$\pi^{-1} (Q)$ i.e. $\pi(C)$ corresponds topologically to a closed
interval of $C'$ and $Q$ is an end-point of this interval. 

Conversely, let $C$
be a connected component of $X(\RR)$ such that $\pi (C)$ is not a
connected component of $X'(\RR)$. Then 
$\pi(C)$, corresponds topologically to a closed
interval of a connected component of $X'(\RR)$. Let $Q$ be
one of the two end-points of this interval. Let $P\in \pi^{-1} (Q)\cap
C$. Then $e_P$ is even since $C$ is clearly on one side of the fiber
$\pi^{-1} (Q)$. Then $P$ is a ramification point of $\pi$.
\end{proof}

\begin{defn}
Let $G$ be a subgroup of $\Aut(X)$. Let $X'$ be the
quotient real curve $X/G$. We define $\Int (\pi )$ as the number of
connected components of the set
\[
\begin{aligned}
  \{\pi(C),&\quad\text{$C$ is a connected component of $X(\RR)$ and}\\
           &\quad\text{$\pi(C)$ is not a connected component of $X'
(\RR)$}\}
\end{aligned}
\]
By the previous proposition, we have  $2\Int (\pi )\leq w_{\RR}
(\pi )$. Let $C$ be a connected component of $X(\RR)$. We denote by
$\Stab (C)$ the stabilizer subgroup of the component $C$ in $G$
i.e. the set of $f\in G$ such that $f(C)=C$.
\end{defn}

From Proposition \ref{basic1}, we derive the following consequence
concerning the real fixed points of an automorphism.

\begin{lem}
\label{basic2}
Let $\varphi$ be an automorphism of $X$ of order $N$. The real
ramification points of $\pi: X\rightarrow X/\langle\varphi\rangle$ 
are in 1 to 1
correspondence with the real fixed
points of $\varphi^{\frac{N}{2}}$ and $\mu_{\RR} (\varphi^i )=0$ for any $i\in\{
1,\ldots,N-1\}$ such that $i\not= \frac{N}{2}$. 
Moreover, if a real ramification point of
$\pi: X\rightarrow
X/\langle\varphi\rangle$ belongs to a connected component $C$ of
$X(\RR )$ then exactly $2$ real ramification points of
$\pi: X\rightarrow
X/\langle\varphi\rangle$ belong to $C$.
\end{lem}

\begin{proof} Let $P$ be a real
ramification point of $\pi: X\rightarrow X/\langle\varphi\rangle$. Then $e_P =2$ by
Proposition \ref{basic1} (\cal{iii}). Consequently
$\Stab(P)$ is a subgroup of order $2$ of $\langle\varphi\rangle$ i.e. $P$ is a real fixed
point of $\varphi^{\frac{N}{2}}$. Since $\Stab(P)=\langle\varphi^{\frac{N}{2}}\rangle$, it
follows that $\mu_{\RR} (\varphi^i )=0$ for any $i\in\{
1,\ldots,N-1\}$ such that $i\not= \frac{N}{2}$. If $P$ is a real fixed
point of $\varphi^{\frac{N}{2}}$, then the image of the connected component
$C$ of $X(\RR )$ containing $P$ by the morphism of degree $2$ 
$\pi':X\rightarrow X/\langle\varphi^{\frac{N}{2}}\rangle$ is a closed interval of a
connected component of $(X/\langle\varphi^{\frac{N}{2}}\rangle) (\RR)$ and $P$
belongs to a fiber
above one end-point of $\pi' (C)$ (see Proposition
\ref{basic1}). Clearly $C$ contains exactly $2$ real
ramification points of $\pi: X\rightarrow X/\langle\varphi\rangle$; $P$ and the point
in the fiber above the other end-point of $\pi' (C)$.
\end{proof}

We give now a generalisation of \cite[Th. 2.2.4]{Kn}.
\begin{prop}
\label{knight}
Let $G$ be a subgroup of $\Aut(X)$. We denote by $X'$ the
quotient real curve $X/G$. Let $C$ be a connected component of $X(\RR)$.
\begin{description}
\item[({\cal{i}})] If $\pi (C)$ is a connected component of $X'(\RR )$
  then $\Stab (C)$ is cyclic.
\item[({\cal{ii}})] If $\pi (C)$ is not a connected component of $X'(\RR
  )$ then $\Stab (C)$ is one of the following groups $\ZZ /2$,
  $\ZZ /2 \times\ZZ /2 $, the diedral group $D_n$ of order $2n$, $n\geq 3$. 
\item[({\cal{iii}})] If $\pi (C)$ is not a connected component of
  $X'(\RR)$ and if $|\Stab (C)|=2n$ then $C$ contains exactly $2n$ real
  ramification points of $\pi$ equally shared between two
  fibers. Moreover, if $P\in C$ is a real fixed point of $\varphi$, then
  the other real fixed point of $\varphi$ in $C$ lies in the same
  fiber as $P$ if and only if $n$ is even.
\end{description}
\end{prop}

\begin{proof}
From \cite[Th. 2.2.4]{Kn}, we know that $\Stab (C)$ is either cyclic or
diedral (only cyclic if $\pi (C)$ is a connected component of $X'(\RR
)$). So we only have to prove statement (\cal{iii}). 

We consider the case $\pi (C)$ is not a connected component of
$X'(\RR)$. 
If $\Stab (C)\simeq \ZZ/2$ the result follows from Proposition
\ref{basic1} and Lemma \ref{basic2}.
Assume $\Stab (C)$ is the diedral group
$D_n =\langle\sigma ,\,\rho \rangle$, $n\geq 2$, with $\sigma$ and $\rho$
corresponding
respectively to a symmetry and
a rotation of order $n$ of the regular polygon with $n$ edges. 
If $(a,b,c,d)\in (\ZZ /n)^4$, we
have $(\rho^a \sigma^b )(\rho^c \sigma^d)=\rho^{a+c(-1)^b }\sigma^{b+d}$.
Since $\pi (C)$ is an interval of $X'(\RR)$, there is a real branch
point $Q$ of $\pi$ such that $C$ contains a real fixed point $P$ of
an element of $\Stab (C)$. By Proposition \ref{basic1}, we have $e_P
=2$. 

Firstly, we assume that $P$ is a fixed point of a symmetry and
without loss of generality we can assume that $\sigma(P)=P$. Hence
$\pi^{-1} (Q)\cap C=\{\rho^k (P),\,\,k=0,\ldots,n-1\}$ and $\rho^k (P)$,
$k\in\{0,\ldots,n-1\}$, is
a fixed point of $\rho^k \sigma\rho^{-k}$. 
If $n$ is even, the other fixed point of $\sigma$
in $C$, 
$\rho^{\frac{n}{2}} (P)$, is contained in $\pi^{-1} (Q)\cap C$ and
conversely.

Secondly, we assume that $n$ is even and that $\rho^{\frac{n}{2}}
(P)=P$. For any $k\in\{0,\ldots,n-1\}$, $\rho^k (P)$ is a fixed point
of $\rho^k \rho^{\frac{n}{2}}\rho^{-k}=\rho^{\frac{n}{2}}$ and $\rho^k
(P)\in C$, since $\rho^k \in\Stab (C)$. From Lemma \ref{basic2}, we
conclude that either $n=2$ or $n=4$. If $n=4$, the other fixed point $\rho
(P)$ of
$\rho^2$ in $C$ is contained in $\pi^{-1} (Q)\cap C$. If $n=2$
the other fixed point $\sigma
(P)$ of
$\rho$ in $C$ is contained in $\pi^{-1} (Q)\cap C$.
\end{proof}

\subsection{Real fixed points of automorphisms of real curves}
We will now study morphisms of degree $2$ between M-curves.
\begin{prop}
\label{degre2} 
Let $\varphi$ be a non-trivial automorphism of order $2$ 
of an M-curve such that 
$\mu_{\RR} (\varphi)>0$. Let $X'$ denote the quotient curve
$X/\langle\varphi\rangle$. then
\begin{description}
\item[({\cal{i}})] $X'$ is an $M-$curve;
\item[({\cal{ii}})] all the fixed points
  of $\varphi_{\CC}$ are real and $\mu (\varphi_{\CC} )=\mu_{\RR}
  (\varphi)=2g+2-4g'$
  where $g'$ denotes the genus of $X'$; 
\item[({\cal{iii}})] all the branch points of
  $\pi_{\CC}$ are real and $w_{\RR} (\pi )=2 \Int
  (\pi)=2g+2-4g'$. Moreover, all these branch points are
  contained in a unique connected component of $X'(\RR)$. Concerning
  the other connected components of $X' (\RR)$, the inverse image by
  $\pi$ of
  each of these is a disjoint union of $2$
  connected components of $X (\RR)$.
\item[({\cal{iv}})] $\sigma (\varphi)$ may be written as a product of
  $g'$ disjoint transpositions.
\end{description}
\end{prop}

\begin{proof} From (\ref{equ3}), we get $\mu (\varphi_{\CC} )=2g+2-4g'$.
Since $|\varphi|=2$, we have $w_{\RR} (\pi )+2w_{\CC} (\pi )=w(\pi_{\CC}
)=2g+2-4g'$ and $w_{\RR} (\pi )=2\Int (\pi )$.
It means that there are $\Int (\pi )$ connected components of $X( \RR
)$ such that the image by $\pi$ of each of these components is not a
connected component of $X' (\RR )$ and $\Int (\pi )\leq g+1-2g'$.
Consequently, there are $g+1-\Int (\pi )$
connected components of $X( \RR
)$ such that the image by $\pi$ of each of these components is a
connected component of $X' (\RR )$ and $g+1-\Int (\pi)\geq 2g'$. 
Let $s'$ denote the number of
connected components of $X' (\RR)$.
Let $v$ denote the number of
connected components of $X' (\RR)$ 
that do not contain any real branch point of $\pi$. According to the
above remarks, we have $v\geq \frac{g+1-\Int (\pi)}{2}\geq
g'$. 
Since $\pi$ has at least
$2$ real branch points and by Harnack inequality, we get $v=g'$,
$s'=v+1$ and $\Int(\pi )=g+1-2g'$. Hence $w_{\RR} (\pi )=w(\pi_{\CC}
)=2g+2-4g'$ and the
statements ({\cal{i}}), ({\cal{ii}}) and
({\cal{iii}}) 
follow. Statement ({\cal{iv}) is a consequence of 
statement  (\cal{iii}}).
\end{proof}

Let us mention two nice consequences of Proposition \ref{degre2}.

\begin{thm}
\label{fixm}
Let $\varphi$ be a non-trivial automorphism of an M-curve. 
If one of the fixed points of $\varphi_{\CC }$ is real
then all are real.
\end{thm}

\begin{proof}  If $\varphi$ has a real
fixed point then $|\varphi|=2$ by Lemma \ref{basic2}. By Proposition
\ref{degre2}, the proof is done.
\end{proof} 

\begin{cor}
\label{fixfree}
Let $\varphi$ be an automorphism of order $N>2$ of an M-curve. If $\pi
:X\rightarrow X/\langle\varphi\rangle$ has at least one real ramification point
then $\varphi_{\CC}$ is fixed point free.
\end{cor}

\begin{proof} Assume $\pi
:X\rightarrow X/\langle\varphi\rangle$ has at least one real
ramification 
point. By Lemma
\ref{basic2}, $N$ is even,
the real ramification points of $\pi:
X\rightarrow X/\langle\varphi\rangle$ are the real fixed points of
$\varphi^{\frac{N}{2}}$ and $\mu_{\RR} (\varphi^i )=0$ for any $i\in\{
1,\ldots,N-1\}$ such that $i\not= \frac{N}{2}$. Let $P$ be a
fixed point of $\varphi_{\CC}$. Since $N>2$, $P$ is not a real point.
Clearly $P$ is also a fixed point of $\varphi_{\CC }^{\frac{N}{2}}$, which
contradicts Theorem \ref{fixm}.
\end{proof}

We state now a generalization of Theorem \ref{fixm}.

\begin{thm}
\label{fixs}
Let $\varphi$ be a non-trivial automorphism of order $N$ of a real
curve. 
If $\pi
:X\rightarrow X/\langle\varphi\rangle$ has at least one real ramification point
then $\mu (\varphi_{\CC})\leq
2+\frac{2g-2s+\mu_{\RR} (\varphi^{\frac{N}{2}}
)}{|\varphi|-1}.$ Consequently, 
$\mu (\varphi_{\CC})-\mu_{\RR} (\varphi^{\frac{N}{2}})\leq
2+\frac{2g-2s}{|\varphi|-1}\leq 2(g+1-s).$ 
If $|\varphi|\geq g+2-s$ then $\mu (\varphi_{\CC})-\mu_{\RR}
(\varphi^{\frac{N}{2}})\leq 2$.
\end{thm}

\begin{proof} By Lemma
\ref{basic2}, $N$ is even,
the real ramification points of $\pi:
X\rightarrow X/\langle\varphi\rangle$ are the real fixed points of
$\varphi^{\frac{N}{2}}$ and $\mu_{\RR} (\varphi^i )=0$ for any $i\in\{
1,\ldots,N-1\}$ such that $i\not= \frac{N}{2}$. 
We denote by $X'$ the quotient space $X/
\langle\varphi\rangle$ and by
$g'$ the genus of $X'$.

By Lemma \ref{basic2}, there are exactly $\frac{\mu_{\RR} (\varphi^{\frac{N}{2}}
)}{2}$ connected components of $X( \RR
)$ containing at least a real ramification point of $\pi$. 
More precisely, each of these connected 
components contains exactly two real ramification points of $\pi$ 
(see Lemma \ref{basic2}).
It means that there are $\frac{\mu_{\RR} (\varphi^{\frac{N}{2}}
)}{2}$ connected components of $X( \RR
)$ such that the image by $\pi$ of each of these components is not a
connected component of $X' (\RR )$. Consequently, there are 
$s-\frac{\mu_{\RR} (\varphi^{\frac{N}{2}})}{2}$
connected components of $X( \RR
)$ such that the image by $\pi$ of each of these components is a
connected component of $X' (\RR )$. Let $v$ denote the number of
connected components of $X' (\RR)$ 
that do not contained any real branch point of $\pi$. By a
previous computation, we have $v\geq \frac{s-\frac{\mu_{\RR} (\varphi^{\frac{N}{2}}
)}{2}}{|\varphi|}$. Since $\pi$ has at least
$2$ real branch points, from Harnack inequality, we get 
$$g'\geq v \geq \frac{s-\frac{\mu_{\RR} (\varphi^{\frac{N}{2}}
)}{2}}{|\varphi|}$$ i.e.
\begin{equation}
\label{equ10} s-\frac{\mu_{\RR} (\varphi^{\frac{N}{2}}
)}{2}\leq g'|\varphi|.
\end{equation}
From (\ref{equ3}), we obtain  $\mu (\varphi_{\CC})\leq
2+\frac{2g-2|\varphi|g'}{|\varphi|-1}$. Combining the previous inequality and 
(\ref{equ10}), we get $$\mu (\varphi_{\CC})\leq
2+\frac{2g-2s+\mu_{\RR} (\varphi^{\frac{N}{2}}
)}{|\varphi|-1}.$$
The rest of the proof follows easily from the previous inequality.
\end{proof}

\subsection{Non-real fixed points of automorphisms of real curves}

The following theorem gives an upper bound on the number of non-real fixed
points of an automorphism in terms of the number of connected
component of the real part of the curve.

\begin{thm}
\label{boundnonreal}
Let $\varphi$ be a non-trivial automorphism of $X$ such that $\pi :X\rightarrow
X'=X/\langle\varphi\rangle$ is without real ramification points.
Then $$\mu (\varphi_{\CC})\leq 4+2\frac{g+1-s}{|\varphi|-1} \,\,{\rm
  if}\,\,s>1$$ and 
$$\mu (\varphi_{\CC})\leq 4+2\frac{g-1}{|\varphi|-1} \,\,{\rm
  if}\,\,s=1.$$
Consequently, $\mu (\varphi_{\CC})\leq 4$
if $|\varphi|\geq g+3-s$ (resp.$|\varphi|\geq g+1$) and if $s>1$ (resp. $s=1$).
\end{thm}

\begin{proof}
Let $s'$ denote the number of connected components of $X'
(\RR)$. Since $\pi :X\rightarrow
X'=X/\langle\varphi\rangle$ is without real ramification points, 
the image of any
connected component of $X(\RR )$ is a connected component of $X'(\RR
)$ (Proposition \ref{basic1} ({\cal{iv}}))) i.e. $\Int (\pi)=0$.
Let $g'$ denote the
genus of $X'$.
By Harnack inequality,
\begin{equation}
\label{equ5}
s\leq |\varphi|s' \leq |\varphi|(g' +1).
\end{equation}
From (\ref{equ3}) and (\ref{equ5}), we obtain  respectively 
$\mu (\varphi_{\CC})\leq
2+\frac{2g-2|\varphi|g'}{|\varphi|-1}$ and $|\varphi|g'\geq
s-|\varphi|$ (since $|\varphi|\geq 2$ we replace the last inequality
by $|\varphi|g'\geq 2-|\varphi|$ in the case $s=1$). 
Combining
the two previous inequalities, we get
\begin{equation}
\label{equ6}
\mu (\varphi_{\CC})\leq 4+2\frac{g+1-s}{|\varphi|-1}
\end{equation}
if $s>1$, and 
\begin{equation}
\label{equ6s=1}
\mu (\varphi_{\CC})\leq 4+2\frac{g-1}{|\varphi|-1}
\end{equation}
if $s=1$.
Since $\varphi$ is real, $\mu (\varphi_{\CC})$ is even. 
By (\ref{equ6}) (resp.(\ref{equ6s=1})), $\mu
(\varphi_{\CC})\leq 4$ if $\frac{g+1-s}{|\varphi|-1}<1$ and $s>1$
(resp. $\frac{g-1}{|\varphi|-1}<1$ and $s=1$)
i.e. if $|\varphi|\geq g+3-s$ (resp. $|\varphi|\geq g+1$).
\end{proof}

Let us state a nice concequence of the previous theorem for
M-curves.
\begin{cor}
\label{boundnonrealmax}
Let $\varphi$ be a non-trivial automorphism of an M-curve $X$ 
such that $\pi :X\rightarrow
X'=X/\langle\varphi\rangle$ is without real ramification points. 
Then $\mu (\varphi_{\CC})\leq 4$.
\end{cor}

For an automorphism $\varphi$ of order $|\varphi |\geq g+2$ 
we may improve the result
of Theorem \ref{boundnonreal}.
\begin{prop}
\label{boundnonrealone}
Let $\varphi$ be an automorphism of a curve $X$ such that
$\pi :X\rightarrow
X'=X/\langle\varphi\rangle$ is without real ramification points.
If $|\varphi|\geq g+2$ then $\mu (\varphi_{\CC})\leq 2$.
\end{prop}

\begin{proof}
From (\ref{equ3}),
we obtain $\mu (\varphi_{\CC})\leq
2+\frac{2g-2|\varphi|g'}{|\varphi|-1}\leq
2+\frac{2g}{|\varphi|-1}$. If $|\varphi|\geq g+2$ then we get $\mu
(\varphi_{\CC})\leq 2$.
\end{proof}

\section{An upper bound on the order of some automorphisms groups of
  real curves with real ramification points}

We give an upper bound for the order of an automorphism $\varphi$ such that
$\pi: X\rightarrow X/\langle\varphi\rangle$ has at least one real ramification point.

\subsection{The cyclic case}

\begin{thm}
\label{ordre1}
Let $\varphi$ be an automorphism of $X$ of order $N>1$ such that \\
$\pi:
X\rightarrow X/\langle\varphi\rangle$ has at least one real ramification point.
Then $$N=\frac{\mu_{\RR} (\varphi^{\frac{N}{2}})}{\Int (\pi)} 
\leq \inf\left\{ \frac{2s}{\Int (\pi)}, \frac{2g+2-4g''}{\Int (\pi)}
\right\}\leq 2g+2$$ where $g''$ denotes the
genus of $X/\langle\varphi^{\frac{N}{2}}\rangle$. 
The decomposition of $\sigma (\varphi)$ contains
$\Int (\pi)$ disjoint cycles of order
$\frac{N}{2}$.
If a connected component $C$ of $X(\RR)$ contains a
real ramification point of $\pi$ then $\Stab (C)=\ZZ /2\ZZ$.
\end{thm}

\begin{proof}
By Lemma
\ref{basic2}, the real ramification points of $\pi:
X\rightarrow X/\langle\varphi\rangle$ are the real fixed points of 
$\varphi^{\frac{N}{2}}$
and $\mu_{\RR} (\varphi^i )=0$ for any $i\in\{
1,\ldots,N-1\}$ such that $i\not= \frac{N}{2}$. Let $r$ denote the number of
real ramification points of $\pi$. According to the above remarks, we have
$r=\mu_{\RR} (\varphi^{\frac{N}{2}})$.
The real ramification points of $\pi$ are contained
in fibers above $2\Int (\pi )$
real branch points of $X/\langle\varphi\rangle$ which correspond to
end-points of the
image by $\pi$ of some connected components of $X(\RR )$ (see
Proposition \ref{basic1}). In each fiber
above one of these real branch points of $X/\langle\varphi\rangle$, there are
$\frac{N}{2}$ real ramification points of $\pi$ contained in distinct
connected components of $X (\RR)$ (see Lemma \ref{basic2}) and it is easy to
check that $\varphi$ operates on
these $\frac{N}{2}$ real ramification points as a cycle of order
$\frac{N}{2}$. 
Hence
$$r=2\Int (\pi )\frac{N}{2},$$ and
the decomposition of $\sigma (\varphi)$ contains
$\Int (\pi)$ disjoint cycles of order
$\frac{N}{2}$.
By Proposition \ref{basic2}, $r\leq 2s$ hence $N\leq \frac{2s}{\Int
  (\pi)}\leq 2s $ since $\Int (\pi )\geq 1$. 
By (\ref{equ3}), $\mu (\varphi_{\CC}^{\frac{N}{2}})\leq 2g+2-4g''$.
Hence $N=\frac{r}{\Int (\pi)}\leq \frac{2g+2-4g''}{\Int (\pi)}$.
By Proposition \ref{knight}, if a connected component $C$ of 
$X(\RR)$ contains a
real ramification point of $\pi$ then $\Stab (C)=\ZZ /2\ZZ$.
\end{proof}

We will now look at the case of an automorphism of maximum order.
\begin{thm} 
\label{maxramreel}
Let $\varphi$ be an automorphism of $X$ of order $2g+2$ such that $\pi:
X\rightarrow X/\langle\varphi\rangle$ has at least one real ramification point. Then $X$
is an hyperelliptic M-curve of even genus, $X/\langle\varphi\rangle\simeq
\PP_{\RR}^1$, $\Int (\pi)=1$, $w_{\RR} (\pi
)=2$ and the ramification index of the fibers over the $2$ real branch
points is $2$, $w_{\CC} (\pi)=1$ and the ramification index of the fibers
over the $2$ conjugate non-real branch points of $\pi_{\CC}$ is $g+1$.
Moreover $\sigma (\varphi)$ is the
cyclic permutation $(1\,2 \ldots g+1)$, $\varphi^{g+1}$ is the
hyperelliptic involution and $\Stab (C)=\ZZ
/2\ZZ=\langle\varphi^{g+1}\rangle$ for any connected component $C$ of $X(\RR)$.
\end{thm}

\begin{proof} Since $|\varphi|=2g+2$, Theorem \ref{ordre1} implies $s=g+1$
and $X/\langle\varphi^{\frac{|\varphi|}{2}}\rangle\simeq
\PP_{\RR}^1$. Consequently, $X$
is an hyperelliptic M-curve and $\varphi^{g+1}$ is the hyperelliptic involution.
Moreover, $\Int (\pi)=\frac{\mu_{\RR} (\varphi^{g+1})}{2g+2}=1$, 
$\sigma (\varphi)$ is the
cyclic permutation $(1\,2 \ldots g+1)$ and $\Stab (C)=\ZZ
/2\ZZ=\langle\varphi^{g+1}\rangle$ for any connected component $C$ of
$X(\RR)$, again by
Theorem \ref{ordre1}. Since $X/\langle\varphi^{g+1}\rangle\simeq
\PP_{\RR}^1$ then clearly $X/\langle\varphi\rangle\simeq
\PP_{\RR}^1$.

Assume $g$ is odd and consider the morphism of degree $4$,
$\pi''':X\rightarrow X/\langle\varphi^{\frac{g+1}{2}}\rangle=X'''$. 
Let $g'''$ denote
the genus of $X'''$. By Theorem \ref{fixm} and Lemma \ref{basic2},
$\varphi_{\CC}^{\frac{g+1}{2}}$ is fixed point free and $\mu
(\varphi_{\CC}^{g+1} )-\mu_{\RR} (\varphi^{g+1} )=0$.
Consequently, the ramification
points of $\pi'''$ are the $2g+2$ real fixed points of
$\varphi^{g+1}$. According to Theorem \ref{ordre1} and 
Proposition \ref{degre2}, we
have
$w_{\RR}(\pi''')=2\Int (\pi''')=g+1$.
Writting (\ref{equ3}) for $\varphi^{\frac{g+1}{2}}$, we have
$2g-2=4(2g'''-2)+2g+2$ i.e. $1=2g'''$, which is impossible. Hence $g$ is even.

Using the results of \cite{Kr-Ne} and since 
$X/\langle\varphi\rangle\simeq \PP_{\RR}^1$,
we see that $\pi$ does not have any real branch point such that the
fiber over this point is non-real. Consequently $w_{\RR} (\pi )=2\Int (\pi
)=2$. Let $Q_1 ,\ldots ,Q_{w_{\CC} (\pi)}$ denote the non-real branch
points of $\pi$. The Riemann-Hurwitz relation (\ref{equ2}) gives
$$2g-2=(2g+2)( -2) + 2g+2
+2(2g+2)\sum_{i=1}^{w_{\CC} (\pi)} (1-\frac{1}{e(Q_i)}).$$
Hence 
\begin{equation}
\label{equ4}
\sum_{i=1}^{w_{\CC} (\pi)} (1-\frac{1}{e(Q_i)})=\frac{g}{g+1}<1.
\end{equation}
If $e(Q_i)=2$ for a $i\in \{ 1,\ldots ,w_{\CC} (\pi)\}$ then, 
writting $Q_i =Q'_i +\bar{Q'_i}$ on $X_{\CC }$,
the points $Q'_i ,\bar{Q'_i}$ are fixed points of $\varphi_{\CC}^{g+1}$ in
contradiction with Theorem \ref{fixm}. Hence $e(Q_i)>2$ for
$i=1,\ldots ,w_{\CC} (\pi)$. By (\ref{equ4})  we conclude that $w_{\CC}
(\pi)=1$, hence that
$e(Q_1)=g+1$. Since, if $g$
is odd, a fixed point of
$\varphi_{\CC}^{2}$ is also a fixed point of $\varphi_{\CC}^{g+1}$, 
it also follows from Theorem \ref{fixm} that $g$ is even.
\end{proof}

We give now some remarks concerning the genus of
$X/\langle\varphi^{\frac{|\varphi|}{2}}\rangle$.

\begin{prop} 
\label{remgenre}
Let $\varphi$ be an automorphism of $X$ of order $N>1$ such that $\pi:
X\rightarrow X/\langle\varphi\rangle$ 
has at least one real ramification point. Let
$g''$ denote the genus of $X/\langle\varphi^{\frac{N}{2}}\rangle$.
If $X/\langle\varphi\rangle\simeq \PP_{\RR}^1$ then $g''\leq \frac{g+1-s}{2}$.
\end{prop}

\begin{proof} Since $X/\langle\varphi\rangle\simeq \PP_{\RR}^1$ and since $\pi$ 
has at least one real ramification point, it follows from Proposition
\ref{basic1}
that the image by $\pi$ of
any connected component of $X(\RR)$ is strictly contained in
$\PP_{\RR}^1 (\RR)$. By Lemma \ref{basic2}, 
any connected component of $X(\RR )$ contain exactly $2$
real ramification points of $\pi$. Since the real ramification
points of 
$\pi:
X\rightarrow X/\langle\varphi\rangle$ are the real fixed points of
$\varphi^{\frac{N}{2}}$ (Lemma \ref{basic2}), we conclude that $2s\leq \mu
(\varphi_{\CC}^{\frac{N}{2}})=2g+2-4g''$, which proves the proposition.
\end{proof}

For a real curve $X$ with many real components, the previous proposition
yields information about the automorphisms $\varphi$ of $X$ such that 
$X/\langle\varphi\rangle\simeq \PP_{\RR}^1$.
\begin{cor}
Let $\varphi$ be an automorphism of $X$ of order $N$ such that $\pi:
X\rightarrow X/\langle\varphi\rangle$ has at least one real ramification point.
If $X/\langle\varphi\rangle\simeq \PP_{\RR}^1$ and $X$ has many real 
components then $X$ is hyperelliptic.
\end{cor}

For an automorphism of prime order, the existence of a real ramification point
for the quotient map forces this automorphism to be an involution.
\begin{prop}
\label{primer}
Let $\varphi$ be an automorphism of $X$ of prime order $p$ such that $\pi:
X\rightarrow X/\langle\varphi\rangle$ has at least one real
ramification point. Then
$p=2$.
\end{prop}

\begin{proof}
If $\pi:
X\rightarrow X/\langle\varphi\rangle$ has at least 
one real ramification point then $\langle\varphi\rangle$
has a subgroup of order $2$. Hence $p=2$.
\end{proof}

\subsection{The abelian case}

We give an upper bound for the order of an abelian group of
automorphisms $G$ such that
$\pi: X\rightarrow X/G$ has at least one real ramification point.

\begin{thm}
\label{boundabreal}
Let $G$ be an abelian group of automorphisms of $X$ of order $N>1$
such that 
$\pi:
X\rightarrow X/G$ has at least one real ramification point.
Then 
$$|G|\leq \inf\left\{ \frac{4s}{\Int (\pi)}, 2g+2+4(g+1-s) \right\}\leq 3g+3.$$
\end{thm}

\begin{proof}
Let $g'$ denote the genus of $X'$.
The real ramification points of $\pi:
X\rightarrow X/G$
are contained
in fibers above $2\Int (\pi )$
real branch points of $X/G$ which correspond to end-points of the
image by $\pi$ of some connected components of $X(\RR )$ (see
Proposition \ref{basic1}). Let $P$ be a real ramification point of
$\pi$. By \cite[Lem. 1.1]{Pa}, the stabilizer subgroup of $P$ in
$G$ is cyclic. Let $\varphi$ be a generator of $\Stab(P)$,
then $P$ is
a real fixed point of $\varphi$ and moreover we have $|\varphi|=2$ 
(see Proposition \ref{basic1}). In the fiber above $Q=\pi(P)$, we have
$\frac{N}{2}$ real fixed points of $\varphi$ since points in the same
fiber have conjugate stabilizers and since $G$ is abelian. 
By Lemma \ref{basic2}, a connected
component of $X(\RR)$ contains at most $2$ real fixed points of $\varphi$.
Hence $\pi^{-1} (Q)$ intersects at least $\frac{N}{4}$ distinct connected
components of $X(\RR )$ i.e. $\frac{N}{4}\leq s$. The same conclusion
can be drawn for any ramified fiber with real points, which proves
that $N\leq\frac{4s}{\Int (\pi)}$.

Before finishing the proof we need to make a remark concerning real
branch points with non-real fiber in the case $g'=0$. 
Let $Q$ be a real branch point such
that $\pi^{-1} (Q)$ is non-real. By \cite[Satz 1]{Kr-Ne}, the
decomposition group of $Q$ in $G$ is the diedral group $D_{e(Q)}$. 
Since $G$ is abelian,
we must have $e(Q)=2$.

Now assume $N> 2g+2$. It follows from the beginning of the proof that
$\Int (\pi )=1$. 
By Proposition \ref{basic1} ({\cal{iv}}), 
there exist exactly two real branch points $Q_1$, $Q_2$ of $\pi$
with real fibers and with $e (Q_1 )=e(Q_2 )=2$. Let $P_1$ (resp. $P_2$) be a
point in the fiber $\pi_{\CC}^{-1} (Q_1)$ (resp. $\pi_{\CC}^{-1} (Q_2)$) . 
Let $\varphi_1$  (resp. $\varphi_2$) be a
generator of $\Stab (P_1)$ (resp. $\Stab (P_2)$). Since $G$ is
abelian, 
all the points in the fiber $\pi_{\CC}^{-1}
(Q_1)$ (resp. $\pi_{\CC}^{-1}
(Q_2)$) are fixed points of $\varphi_{1}$ (resp. $\varphi_{2}$). 
Since $N>2g+2$, it follows from the Riemann-Hurwitz formula (\ref{equ2}) that
$g'=0$ and there exist at least two distinct branch points 
$Q'_1$, $Q'_2$ of $\pi_{\CC}$ with non-real fibers. We may assume that
the points $Q'_1$,
$Q'_2$ are either both real or both non-real and conjugate.
Let $\varphi'_1$  (resp. $\varphi'_2$) be a
generator of the stabilizer subgroup of any point in $\pi_{\CC}^{-1}
(Q'_1)$ (resp. $\pi_{\CC}^{-1}
(Q'_2)$). We have different cases.

{\bf Case 1:} $\varphi_{1}$ and $\varphi_{2}$ are powers of either of
$\varphi'_1$ or $\varphi'_2$.\\
By Theorem \ref{fixs}
and since there are $\frac{|G|}{e(Q'_1)}$ points in
the fiber $\pi_{\CC}^{-1} (Q'_1)$, we obtain
$\frac{|G|}{e(Q'_1)} \leq 2(g+1-s)$ and $X$ is not an M-curve.
Similarly, we have $\frac{|G|}{e(Q'_2)}\leq 2(g+1-s)$.
By (\ref{equ2}) and according to above remarks, 
we obtain $2g-2\geq -2|G|+\sum_{i=1}^2
|G|(1-\frac{1}{e(Q_i)})+\sum_{i=1}^2
(|G|-\frac{|G|}{e(Q'_i)})\geq |G|-4(g+1-s)$. Hence $|G|\leq
2g-2+4(g+1-s)$.

{\bf Case 2:} $\varphi_{1}$ is a power of $\varphi'_1$, 
and $\varphi_2$ is not a power of either of $\varphi'_1$ and $\varphi'_2$.\\
By Theorem \ref{fixs}, we have
$\frac{|G|}{e(Q_1)} \leq \mu(\varphi'_{1,\CC})-\mu_{\RR} (\varphi_1)\leq
2(g+1-s)$ and $X$ is not an M-curve.
Since $\varphi_2$ is not a power of either of $\varphi'_1$ and $\varphi'_2$ and
since the fibers 
$\pi_{\CC}^{-1} (Q_1)$ and $\pi_{\CC}^{-1} (Q_2)$ 
contain all the real fixed points of the elements of $G$, it follows that
the map $X\rightarrow X/\langle\varphi'_2\rangle$ 
has no real ramification point.
By (\ref{equ6}), we have $\frac{|G|}{e(Q'_2)}\leq
\mu(\varphi'_{2,\CC})\leq 4+2(g+1-s)$.
By (\ref{equ2}) and according to the above remarks, 
we obtain $2g-2\geq -2|G|+\sum_{i=1}^2
|G|(1-\frac{1}{e(Q_i)})+\sum_{i=1}^2
(|G|-\frac{|G|}{e(Q'_i)})\geq |G|-4-4(g+1-s)$. Hence $|G|\leq
2g+2+4(g+1-s)$.

{\bf Case 3:} $\varphi_{1}$ and $\varphi_{2}$ are not powers of either of
$\varphi'_1$ or $\varphi'_2$ and $e(Q'_1 )\geq 3$.\\
Similarly to the previous case, the maps 
$X\rightarrow X/\langle\varphi'_1\rangle$  and 
$X\rightarrow X/\langle\varphi'_2\rangle$
are without real ramification points. By an above remark, $Q'_1$ is
non-real and we may assume that $Q'_1$ and $Q'_2$ are two conjugate
points of $X'_{\CC}$ corresponding to the non-real point $Q'$ of $X'$.
By (\ref{equ6}) and since $e(Q' )\geq 3$, 
we have $\frac{|G|}{e(Q')} \leq 2+\frac{g+1-s}{2}$.
By (\ref{equ2}) and according to the above remarks, 
we obtain $2g-2\geq -2|G|+\sum_{i=1}^2
|G|(1-\frac{1}{e(Q_i)})+
2(|G|-\frac{|G|}{e(Q')})\geq |G|-4-(g+1-s)$. Hence $|G|\leq
2g+2+(g+1-s)$. Since we have assumed $N>2g+2$, 
this case does not occur when $X$ is an M-curve.

{\bf Case 4:} $\varphi_{1}$ and $\varphi_{2}$ are not powers of either of
$\varphi'_1$ or $\varphi'_2$ and $e(Q'_1 )=e(Q'_2 )=2$.\\
By (\ref{equ2}), there exists a third branch point $Q'_3$ of
$\pi_{\CC}$ 
with non-real fibers.
By Theorem \ref{fixs} and (\ref{equ6}), we obtain $\frac{|G|}{e(Q'_3)}\leq 4+2(g+1-s)$.
By (\ref{equ2}), 
we get $2g-2\geq -2|G|+\sum_{i=1}^2
|G|(1-\frac{1}{e(Q_i)})+\sum_{i=1}^3
(|G|-\frac{|G|}{e(Q'_i)})\geq |G|-4-2(g+1-s)$. Hence $|G|\leq
2g+2+2(g+1-s)$. Since we have assumed $N>2g+2$, 
this case does not occur when $X$ is an M-curve.

We have proved that $|G|\leq \inf\{ 4s, 2g+2+4(g+1-s) \}$. Since
$4s\geq 2g+2+4(g+1-s)$ if and only if $s\geq \frac{3}{4}(g+1)$, we get
$|G|\leq 3g+3$.
\end{proof}

For M-curves, we obtain:
\begin{cor}
Let $G$ be an abelian group of automorphisms of an M-curve $X$ such that $\pi:
X\rightarrow X/G$ has at least one real ramification point.
Then 
$|G|\leq 2g+2 $.
\end{cor}

\subsection{The hyperelliptic case}

We will give an upper bound for the order of the group of
automorphisms of a real hyperelliptic curve such that the
hyperelliptic involution
has at least one real fixed point. Before that, we will prove a more
general result.

\begin{thm}
\label{boundcenterreal}
Let $G$ be a group of automorphism of $X$ of order $N>1$ such that there
exists $\varphi\not= Id$ in the center $Z(G)$ of $G$ 
with at least one real fixed point.
Then $|G|\leq \inf \{ 4s,4g+4-8g'' \} \leq 4g+4$ where $g''$
denotes the genus of $X/\langle\varphi\rangle$.
\end{thm}

\begin{proof}
Let $P$ be a a real fixed point of $\varphi$, then $P$ is a real
ramification point of $\pi:
X\rightarrow X/G$ with ramification index $e_P =2$ (see Proposition
\ref{basic1}).
Since the points in the fiber above $Q=\pi (P)$
have conjugate stabilizer subgroups and since $\varphi\in Z(G)$,
the fiber
$\pi^{-1} (Q)$ is composed by $\frac{N}{2}$ real points which are real
fixed points of $\varphi$. By Lemma \ref{basic2}, a connected
component of $X(\RR)$ contains at most $2$ real fixed points of $\varphi$.
Hence $\pi^{-1} (Q)$ intersects at least $\frac{N}{4}$ distinct connected
components of $X(\RR )$ i.e. $\frac{N}{4}\leq s$. 
By (\ref{equ3}),
$\mu_{\RR} (\varphi)\leq 2g+2-4g''$. Hence $\frac{N}{2}\leq 2g+2-4g''$
and the proof is done.
\end{proof}

\begin{cor} 
\label{boundhypreal}
Let $X$ be a real hyperelliptic curve such that the
hyperelliptic involution $\imath$ has at least a real fixed point
(e.g. if $s\geq 3$). Then $|\Aut (X)|\leq 4s$.
\end{cor}

\begin{proof}
By \cite[Cor. 3 p. 102]{Fa-K}, the hyperelliptic involution $\imath$
is in the center of $\Aut (X)$. If $s\geq 3$ then $\imath$ has at
least a real fixed point \cite[Prop. 4.3]{Mo}. The rest of the proof
follows from Theorem \ref{boundcenterreal}.
\end{proof}

\begin{rema}
{\rm The order of the automorphism group of a real curve $X$ cannot be
  larger than $12(g-1)$ \cite{Ma1}. In the case $|\Aut
  (X)|=12(g-1)$, the map $\pi:
  X\rightarrow X/\Aut(X)$ has at least a real ramification point
  \cite{Ma1} and it follows from Corollary \ref{boundhypreal} that $X$
  is not hyperelliptic. In the cyclic case, the curves with an
  automorphism of maximum order are hyperelliptic.}
\end{rema}

We extend the result concerning the hyperelliptic curves to real
curves which are $2$-sheeted coverings.
\begin{cor} 
\label{bound2real}
Let $X$ be a real curve such that $X$ is a $2$-sheeted covering of a
curve of genus $g''$. Assume $g>4g'' +1$ and the involution
induced by the $2$-sheeted covering has at least a real fixed point,
then $|\Aut (X)|\leq \inf \{ 4s,4g+4-8g'' \} \leq 4g+4$.
\end{cor}

\begin{proof} 
If $g>4g'' +1$, the involution induced by the $2$-sheeted covering
is in the center of $\Aut (X)$
\cite[Thm. p. 250]{Fa-K}. The proof
follows now from Theorem \ref{boundcenterreal}.
\end{proof}

\section{An upper bound for some groups of automorphisms of real 
curves without real ramification points}

This section is devoted to the study of automorphisms of real curves
without real fixed points.

\subsection{The cyclic case for M-curves}

In this section we give an upper bound on the order of an automorphism 
of an M-curve using Corollary \ref{boundnonrealmax}.

Before giving this bound, we need to give a result concerning the
number of branch points of a map corresponding to a quotient by an
abelian group of automorphisms.

\begin{lem}
\label{pardini}
Let $G$ be an abelian group of automorphisms of a 
smooth projective curve $X$ over
$\CC$. Let $\pi$ denote the map $\pi :X\rightarrow X'=X/G$. Then
$w (\pi )\not= 1$.
\end{lem}

\begin{proof} 
Write $G$ as the direct sum of the cyclic subgroups $G_1 ,\ldots ,G_t$
generated by $\varphi_1 ,\ldots ,\varphi_t$.
Assume $\pi$ has a unique branch point denoted by $Q$. From the
algebra structure of $\pi_* \O_X$ and from the action of $G$ on it,
it is possible to derive a linear equivalence, for $i=1,\ldots
,t$, between $|\varphi_i|D_i$ and $n_iQ$,
with $D_i$ a divisor on $X'$ associated to the dual of $\langle\varphi_i\rangle$
and $n_i$ an integer such that $0\leq n_i <|\varphi_i|$
\cite[Prop. 2.1]{Pa}.
By \cite[Prop. 2.1]{Pa}, the data: $D_i$, $Q$ and the linear
equivalence between $|\varphi_i|D_i$ and $n_i Q$, determines 
uniquely $\pi$ up to
isomorphism. According to the previous remark and since 
$\pi$ is ramified, we can assume that $n_1 \geq 1$ (see
\cite[Example 2.1 ({\cal{ii}})]{Pa}).
But then the condition $n_1 <|\varphi_1|$ says that the linear
equivalence between $|\varphi_1 | D_1$ and $n_1 Q$ is impossible. 
Hence $w (\pi )> 1$.
\end{proof}

In the following theorem, we determine all the automorphisms of M-curves.
\begin{thm}
\label{ordermcurve}
Let $\varphi$ be a non-trivial automorphism of an M-curve $X$ 
such that $\pi :X\rightarrow
X'=X/\langle\varphi\rangle$ is without real ramification points. 
Let $g'$ denote the
genus of $X'$.
Then one of the five possibilities occurs.
\begin{description}
\item[({\cal{i}})] $\mu (\varphi_{\CC} )=4$, $|\varphi|=\frac{g+1}{g'+1}$,
  $X'$ is an M-curve, $\pi_{\CC}$ has $4$ branch points $Q_1
  ,\bar{Q_1},Q_2 ,\bar{Q_2}$ i.e. $w_{\CC} (\pi)=2$ and $w_{\RR} (\pi
  )=0$, $e(Q_1)=e(Q_2)=|\varphi|$.
\item[({\cal{ii}})] $\mu (\varphi_{\CC} )=2$, $|\varphi|=\frac{g}{g'}$ with
  $g'\geq 1$, $X'$ is an M-curve, $\pi_{\CC}$ has $2$ branch points $Q_1
  ,\bar{Q_1}$ i.e. $w_{\CC} (\pi)=1$ and $w_{\RR} (\pi
  )=0$, $e(Q_1)=|\varphi|$.
\item[({\cal{iii}})] $\mu (\varphi_{\CC} )=0$, $|\varphi|=\frac{g+1}{g'}$ with
  $g'\geq 1$,
  $X'$ has many real components, $\pi_{\CC}$ has $2$ branch points $Q_1
  ,\bar{Q_1}$ i.e. $w_{\CC} (\pi)=1$ and $w_{\RR} (\pi
  )=0$, $e(Q_1)=\frac{|\varphi|}{2}$.
\item[({\cal{iv}})] $\mu (\varphi_{\CC} )=0$, $|\varphi|=\frac{g+1}{g'}$ with
  $g'\geq 1$,
  $X'$ is an M-curve, $\pi_{\CC}$ has $2$ real branch points $Q_1
  ,Q_2$ with non-real fibers i.e. $w_{\CC} (\pi)=0$ and $w_{\RR} (\pi
  )=2$, $e(Q_1)=e(Q_2)=\frac{|\varphi|}{2}$. The real branch points $Q_1
  ,Q_2$ are contained in the same connected component of $X'(\RR)$.
\item[({\cal{v}})] $\mu (\varphi_{\CC} )=0$, $|\varphi|=\frac{g-1}{g'-1}$ with
  $g'\geq 2$, $X'$ has many real components, 
  $\pi_{\CC}$ does not have any branch point.
\end{description}
It follows that $|\varphi|\leq g+1$ and $X'$ has many real components.
\end{thm}

\begin{proof}
Let $s'$ denote the number of connected components of $X'(\RR )$. By
Corollary \ref{boundnonrealmax}, $\mu (\varphi_{\CC} )\leq 4$. We will
proceed by looking successively at the cases $\mu (\varphi_{\CC}
)=4,2$ and $0$.

Assume $\mu (\varphi_{\CC} )=4$. By Corollary \ref{boundnonrealmax} and
since a fixed point of $\varphi_{\CC}$ is fixed for any power of
$\varphi_{\CC}$, the Riemann-Hurwitz relation (\ref{equ3}) reads
$2g-2=|\varphi|(2g'-2)+4(|\varphi|-1)$. It gives
$|\varphi|=\frac{g+1}{g'+1}$. By (\ref{equ5}) we have
$s'\geq\frac{g+1}{|\varphi|}
=g'+1$. Harnack inequality says that $X'$ is an M-curve.  
We are in the case of statement ({\cal{i}}).

Assume $\mu (\varphi_{\CC} )=2$. Since a totally
ramified non-real point of $\pi$ is necessarily in a fiber above a
non-real branch point, $\pi_{\CC}$ has $2$ branch points $Q_1
,\bar{Q_1}$ such that $e(Q_1)=|\varphi|$. Let $P_1$ (resp. $ \bar{P_1}$)
denote the totally ramified point over $Q_1$ (resp. $ \bar{Q_1}$).
If $\pi_{\CC}$ does not have more
branch points, then (\ref{equ3}) gives $|\varphi|=\frac{g}{g'}$ with
$g'\geq 1$. By (\ref{equ5}), $X'$ is an M-curve and we get
statement ({\cal{ii}}).
If $\pi_{\CC}$ has at least one more
branch point then set 
\[
j=\inf\{ \,i>1,\,\, \varphi_{\CC}^i \,\,\text{has a fixed
point $P_2$ and $P_1$ is not fixed by $\varphi_{\CC}^i$\}}. 
\]
Since the
fiber containing  $P_2$ contains $j$ fixed points of $\varphi_{\CC}^j$ and
since $P_1$ and $ \bar{P_1}$ are  fixed points of $\varphi_{\CC}^j$, 
Corollary \ref{boundnonrealmax} implies that $j=2$ and that $P_2$ and 
$ \bar{P_2}$ are in contained in a fiber above a real branch point. 
It follows  that $w (\pi_{\CC})=3$. By (\ref{equ2}), we get $2g-2=
|\varphi|(2g'-2)+2(|\varphi|-1)+|\varphi|-2$ i.e. 
$|\varphi|=\frac{2g+2}{2g'+1}$. 
Since $\pi$ has a real branch point and since $\pi :X\rightarrow
X'=X/\langle\varphi\rangle$ is without real ramification points i.e. 
$\Int (\pi)=0$, 
the real branch point is
contained in a connected component $C'$ of $X'(\RR)$ such that
$\pi^{-1} (C')$ is totally non real. Hence $s'\geq 2$ and $g'\geq 1$
by Harnack inequality.
So we may refine (\ref{equ5}) in
this case and we get $g+1\leq
|\varphi|g'=g'\frac{2g+2}{2g'+1}$. It gives a contradiction
and the single case of an automorphism $\varphi$ with $\mu
(\varphi_{\CC} )=2$ is
given by statement ({\cal{ii}}).

Finally, assume $\varphi_{\CC}$ is fixed point free and 
$\pi_{\CC}$ does not have any branch point.
From (\ref{equ2}) and (\ref{equ5}), it follows that 
$|\varphi|=\frac{g-1}{g'-1}$ with
$g'\geq 2$ and that $X'$ has many real components.
Assume $\varphi_{\CC}$ is fixed point free and $w (\pi_{\CC})\geq 1$. Let
\[ j=\inf\{ \,i>1,\,\, \text{$\varphi_{\CC}^i$ has a fixed
point \}}.\] 
Let $P_1$ be a fixed point of $\varphi_{\CC}^j$. 
Then $ \bar{P_1}$ is also a  fixed point of
$\varphi_{\CC}^j$. Firstly, assume that $P_1$ and $ \bar{P_1}$ are contained
in two conjugate fibers above two conjugate branch points denoted by 
$Q_1$ and $ \bar{Q_1}$. Since the
fiber containing  $P_1$ contains $j$ fixed points of $\varphi_{\CC}^j$, 
Corollary \ref{boundnonrealmax} implies that $j=2$.
By Corollary
\ref{boundnonrealmax}, we see that a ramified fiber of $\pi_{\CC}$
contains
at most $4$ points and the number of points in the fiber is exactly
the smallest power of $\varphi$ which generates the stabilizer 
group of any point in
the fiber. Since a non-real fiber over a real point and two conjugate
fibers contain both an even number of points, we conclude that a
ramified fiber contains $2$ or $4$ points.
Hence $w
(\pi_{\CC})=2$ since if we assumed $w(\pi_{\CC})>2$, we would have $\mu
(\varphi_{\CC}^2 )>4$ or $\mu
(\varphi_{\CC}^4 )>4$ which contradicts Corollary
\ref{boundnonrealmax} (a fixed point of $\varphi_{\CC}^2$ is also
a fixed point of $\varphi_{\CC}^4$).
By (\ref{equ2}), we get $2g-2=
|\varphi|(2g'-2)+2(|\varphi|-2)$ i.e. $|\varphi|=\frac{g+1}{g'}$ and
$g'\geq 1$. By (\ref{equ5}), $X'$ has many real
components and we get statement ({\cal{iii}}).
Assume now that $P_1$ and $ \bar{P_1}$ are contained in the same fiber
above a real branch point denoted by $Q_1$. By Lemma \ref{pardini},
$\pi_{\CC}$ has a least one more branch point. From Corollary
\ref{boundnonrealmax} and according to the above remarks, it follows that 
$\pi_{\CC}$ has exactly $2$ ramified fibers over two real branch
points $Q_1$ and $Q_2$, one fiber contains $P_1$ and $ \bar{P_1}$ and the other
contains the points $P_2$ and 
$ \bar{P_2}$. Moreover, $e(Q_1)=e(Q_2)=\frac{|\varphi|}{2}$. From
(\ref{equ2}), we get $|\varphi|=\frac{g+1}{g'}$ and $g'\geq 1$. 
Arguing as in the
proof of the case $\mu (\varphi_{\CC} )=2$, we see that the inverse image by
$\pi$ of the connected
components of $X' (\RR)$ containing $Q_1$ and $Q_2$ are totally
non-real.
We may refine (\ref{equ5}) in
this case and we get $g+1\leq
|\varphi|(s'-1)\leq |\varphi| g'=g+1$. Hence $X'$ is an M-curve and $Q_1$ and
$Q_2$ are contained in the same connected component of $X'(\RR)$. 
We are in the case of statement ({\cal{iv}}) and the proof is done.
\end{proof}

\subsection{The cyclic case}

At the beginning of this section 
we give an upper bound for the order of an automorphism
of a real curve and we study the limit cases.

\begin{thm}
\label{order}
Let $\varphi$ be a non-trivial automorphism of a real curve $X$ 
such that $\pi :X\rightarrow
X'=X/\langle\varphi\rangle$ is without real ramification points. Then 
$|\varphi|\leq \sup\{2g+4-s, 2g+2-\frac{2}{3}s \}$ if $s>1$ and
$|\varphi|\leq 2g+2$ if $s=1$.
\end{thm}

\begin{proof}
Let $g'$ denote the genus of $X'$. We assume that $|\varphi|\geq
g+3-s$ and by Theorem \ref{boundnonreal} we have $\mu (\varphi_{\CC} )\leq
4$. We will
proceed by looking successively at the cases $\mu (\varphi_{\CC}
)=4,2$ and $0$.
Let $s'$ denote the number of connected components of $X' (\RR)$.

{\bf Case 1:} $\mu (\varphi_{\CC} )=4$.
The Riemann-Hurwitz relation (\ref{equ3}) gives
$2g-2 \geq |\varphi|(2g'-2)+4(|\varphi|-1)$. It gives
$|\varphi|\leq \frac{g+1}{g'+1}\leq g+1$.

{\bf Case 2:}  $\mu (\varphi_{\CC} )=2$.
By (\ref{equ3}), we get $|\varphi|\leq
\frac{g}{g'}\leq g$ if $g'\geq 1$. So assume $g' =0$.
Since $\mu (\varphi_{\CC} )=2$, $\pi_{\CC}$ has already $2$ branch points $Q_1
,\bar{Q_1}$ (i.e. $w_{\CC}(\pi )\geq 1$) 
such that $e(Q_1)=|\varphi|$ (since $s'=1$, the two branch points cannot be
real). Let $P_1$ (resp. $ \bar{P_1}$)
denote the totally ramified point over $Q_1$ (resp. $ \bar{Q_1}$).
From (\ref{equ3}), we see that $\pi_{\CC}$ has at least one more
branch point. Since $\pi :X\rightarrow
X'=X/\langle\varphi\rangle$ is without real ramification points and
since 
$s'=1$,
then $w_{\CC}(\pi )\geq 2$. Let $Q_2$, $ \bar{Q_2}$ denote two conjugate
branch points of $\pi_{\CC}$ distinct from $Q_1$ and $
\bar{Q_1}$. Let $\frac{|\varphi|}{j}=e(Q_2 )$. We have $2\leq
j\leq\frac{|\varphi|}{2}$. The Riemann-Hurwitz relation
(\ref{equ2})
gives $2g-2 \geq -2|\varphi| +2|\varphi|(1-\frac{1}{|\varphi|})
+2|\varphi|(1-\frac{j}{|\varphi|})$
i.e. 
\begin{equation}
\label{equint1}
|\varphi|\leq g+j.
\end{equation} 
The stabilizer subgroups of the points in
the fibers $\pi_{\CC}^{-1} (Q_2 )$ and $\pi_{\CC}^{-1} (\bar{Q}_2 )$
are generated by $\varphi_{\CC}^j$ and each fiber contains $j$
points. Moreover, $P_1$ and $\bar{P_1}$ are also fixed points of
$\varphi_{\CC}^j$ since they are fixed points of
$\varphi_{\CC}$. Consequently, we get $2j+2\leq \mu
(\varphi_{\CC}^j)$. If $s>1$, using (\ref{equ6}) we obtain 
$2j+2\leq 4+2\frac{g+1-s}{\frac{|\varphi|}{j}-1}$. Since $2\leq
\frac{|\varphi|}{j}$, we get $j\leq 1+(g+1-s)$. By (\ref{equint1}), we
get $|\varphi |\leq 2g+2-s$ in the case $s>1$. If $s=1$, 
using (\ref{equ6s=1}) we obtain 
$2j+2\leq 4+2\frac{g-1}{|\frac{\varphi}{j}|-1}$. Similarly to the case
$s>1$, we get $|\varphi |\leq 2g$ in the case $s=1$.

{\bf Case 3:}  $\mu (\varphi_{\CC} )=0$.
Assume $\varphi_{\CC}$ is fixed point free and 
$\pi_{\CC}$ does not have any branch point.
From (\ref{equ2}) and (\ref{equ5}), it follows that 
$|\varphi|=\frac{g-1}{g'-1}\leq g-1$ with
$g'\geq 2$.

In the rest of the proof we assume $\varphi_{\CC}$ is fixed point free and
that
$w (\pi_{\CC})\geq
1$. 

Firstly, we assume that $w_{\CC} (\pi )\geq 1$ i.e. that
$\pi_{\CC}$ has two conjugate
branch points, denoted by $Q_1$, $ \bar{Q_1}$, such that
$\frac{|\varphi|}{j}=e(Q_1 )$ with $2\leq
j\leq\frac{|\varphi|}{2}$. From (\ref{equ2}), we obtain 
$2g-2 \geq |\varphi|(2g'-2) +2|\varphi|(1-\frac{j}{|\varphi|})$
i.e. 
\begin{equation}
\label{equint2}
|\varphi|g'\leq g-1+j.
\end{equation}
The $2j$ points in
the fibers $\pi_{\CC}^{-1} (Q_1 )$ and $\pi_{\CC}^{-1} (\bar{Q}_1 )$
are fixed by $\varphi_{\CC}^j$. Consequently, we get $2j\leq \mu
(\varphi_{\CC}^j)$. If $s>1$ (resp. $s=1$), using (\ref{equ6})
(resp. (\ref{equ6s=1}))
we obtain 
$2j\leq 4+2\frac{g+1-s}{\frac{|\varphi|}{j}-1}\leq 4+2(g+1-s)$
(resp. $2j\leq 4+2\frac{g-1}{\frac{|\varphi|}{j}-1}\leq 4+2(g-1)$). 
If $g'\geq 1$, it follows from (\ref{equint2}) that $|\varphi|\leq
2g+2-s$ if $s>1$ (resp. $|\varphi|\leq 2g$ if $s=1$).

So assume
$g'=0$. By (\ref{equ2}), we see that $\pi_{\CC}$ has at least one more
branch point. Since $\pi :X\rightarrow
X'=X/\langle\varphi\rangle$ is without real ramification points and 
since $s'=1$,
then $w_{\CC}(\pi )\geq 2$. Let $Q_2$, $ \bar{Q_2}$ denote two conjugate
branch points of $\pi_{\CC}$ distinct from $Q_1$ and $
\bar{Q_1}$. Let $\frac{|\varphi|}{j'}=e(Q_2 )$. We have $2\leq
j'\leq\frac{|\varphi|}{2}$. 

If $j=\frac{|\varphi|}{2}$ or 
$j'=\frac{|\varphi|}{2}$, the points in the fibers $\pi_{\CC}^{-1}
(Q_1 )$ and $\pi_{\CC}^{-1} (\bar{Q}_1 )$, or in the fibers $\pi_{\CC}^{-1}
(Q_2 )$ and $\pi_{\CC}^{-1} (\bar{Q}_2 )$, are fixed points of
$\varphi_{\CC}^{\frac{|\varphi|}{2}}$. Hence, using (\ref{equ6})
(resp. (\ref{equ6s=1})) we get $|\varphi|\leq 4+2(g+1-s)=2g+6-2s\leq 2g+4-s$
if $s>1$ (resp. $|\varphi|\leq 2g+2$ if $s=1$).

We assume that $j\leq \frac{|\varphi|}{3}$ and  $j'\leq \frac{|\varphi|}{3}$.
By (\ref{equ2}), we get 
\begin{equation}
\label{equ7}
|\varphi|\leq g-1+j+j'.
\end{equation}
Similarly to the previous cases, we get $2j\leq \mu
(\varphi_{\CC}^j)$ and $2j'\leq \mu
(\varphi_{\CC}^{j'})$. Since $3\leq \frac{|\varphi|}{j}$, it follows
from (\ref{equ6})
(resp. (\ref{equ6s=1})) that $2j\leq 4+(g+1-s)$ if $s>1$
(resp. $2j\leq 4+(g-1)$ if $s=1$) and the same is
true for $j'$. Hence $j+j'<4+(g+1-s)$ if $s>1$ (resp. $j+j'<g+3$ if
$s=1$) 
(equality is impossible in the two previous inequality
since it would imply $j=j' =\frac{|\varphi|}{3}$ and we would get a
contradiction with (\ref{equ6}) and  (\ref{equ6s=1})).
By (\ref{equ7}) and the previous results we get $|\varphi|<
2g+4-s$ if $s>1$ (resp. $|\varphi|< 2g+2$ if $s=1$).

Secondly, we turn to the case $w_{\CC} (\pi )= 0$ and $w_{\RR} (\pi )\geq 1$.
From Lemma \ref{pardini} and since $\pi :X\rightarrow
X'=X/\langle\varphi\rangle$ is without real ramification points, 
we have $w_{\RR} (\pi )\geq
2$ and $s'\geq 2$. In particular, $g'\geq 1$ since $s'\geq 2$.
Let $Q_1$, $Q_2$ be two real branch points of $\pi$. We set 
$\frac{|\varphi|}{j}=e(Q_1 )$ and 
$\frac{|\varphi|}{j'}=e(Q_2 )$. Since the fibers of $\pi$ 
above $Q_1$ and $Q_2$ are non-real, we have $2\leq
j\leq\frac{|\varphi|}{2}$ and $2\leq
j'\leq\frac{|\varphi|}{2}$. By (\ref{equ2}) and taking account of 
the previous remarks, 
we obtain 
\begin{equation}
\label{equint3}
2g' |\varphi|\leq 2g-2+j+j'.
\end{equation}

If $j=\frac{|\varphi|}{2}$ and
$j'=\frac{|\varphi|}{2}$, the $j+j'$ points in the fibers 
$\pi_{\CC}^{-1}
(Q_1 )$ and $\pi_{\CC}^{-1}
(Q_2 )$ are fixed points of
$\varphi_{\CC}^{\frac{|\varphi|}{2}}$. Hence, using (\ref{equ6})
(resp. (\ref{equ6s=1})) we get $j+j'\leq 4+2(g+1-s)$ if $s>1$
(resp. $j+j'\leq 2g+2$ if $s=1$).
Combining (\ref{equint3}) and the previous result, we get
$|\varphi|\leq
2g+2-s$ if $s>1$ (resp. $|\varphi|\leq 2g$ if $s=1$).

If $j\leq\frac{|\varphi|}{2}$ (i.e. $e(Q_1)\geq 2$) and 
$j'<\frac{|\varphi|}{2}$ (i.e. $e(Q_2)\geq 3$). 
The $j'$ points in the fiber
$\pi_{\CC}^{-1} (Q_2)$ are fixed points of 
$\varphi_{\CC}^{\frac{|\varphi|}{j'}}$.
Using (\ref{equint3}), we get
\begin{equation}
\label{equint4}
|\varphi|\leq \frac{4}{3}g-\frac{4}{3}+\frac{2}{3}j'.
\end{equation}
Since $e(Q_2)=\frac{|\varphi|}{j'}\geq 3$, using (\ref{equ6})
(resp. (\ref{equ6s=1})) and  (\ref{equint4}), we obtain
$|\varphi|\leq
2g+2-\frac{2}{3}s$ if
$s>1$ (resp. $|\varphi|\leq 2g+\frac{2}{3}$ if $s=1$). 
\end{proof}

We will now look at the case of an automorphism of maximum order.
\begin{thm} 
\label{maxramnonreel}
Let $\varphi$ be an automorphism of $X$ of order $2g+2$ such that $\pi:
X\rightarrow X/\langle\varphi\rangle$ is 
without real ramification points. Then $X$
is an hyperelliptic curve of even genus, $s=1$, $X/\langle\varphi\rangle\simeq
\PP_{\RR}^1$,  $w_{\RR} (\pi
)=0$ and $w_{\CC} (\pi)=2$. Let $\{ Q_1,Q_2 ,\bar{Q_1},\bar{Q_2}\}$ be the
branch points of $\pi_{\CC}$, then $e(Q_1 )=2$, $e(Q_2 )=g+1$. 
Moreover, $\varphi^{g+1}$ corresponds to the hyperelliptic
involution and $X$ is given by the real polynomial equation $y^2
=f(x)$, where $f$ is a monic polynomial of degree $2g+2$ and where $f$
has no real roots.

\end{thm}

\begin{proof} 
Let $g'$ (resp. $g''$) denote the genus of $X'=X/\langle\varphi\rangle$
(resp. of $X''=X/\langle\varphi^{g+1}\rangle$). Looking at the proof of Theorem
\ref{order}, we see that $g'=0$, $s=1$ or $s=2$ and $\pi_{\CC }$ has exactly
$4$ non-real branch points denoted by $Q_1,Q_2
,\bar{Q_1},  \bar{Q_2}$. If we set $e(Q_1 ) =\frac{|\varphi|}{j}$ and $e(Q_2
)=\frac{|\varphi|}{j'}$, it follows from the proof of Theorem \ref{order}
that we can assume $e(Q_1 )=2$ i.e. $j=g+1$.
By (\ref{equ2}), 
we get $j'=2$
i.e. $e(Q_2 )=g+1$. From the proof of Theorem \ref{order}, we see that
$\mu (\varphi_{\CC}^{g+1} )=|\varphi|=2g+2$. Hence, considering the
relation (\ref{equ3}) 
for $\varphi^{g+1}$, we
get $g''=0$ i.e. $X$ is hyperelliptic and $\varphi^{g+1}$ is the
hyperelliptic involution. If $g$ is odd then the fixed points of 
$\varphi_{\CC}^2$
are fixed points of $\varphi_{\CC}^{g+1}$ and, since $e(Q_1 )=2$ and $e(Q_2
)=g+1$, 
we get $\mu (\varphi_{\CC}^{g+1}
)=2g+6$. It leeds to a contradiction
and we conclude that $g$ is even. Since the map $\pi'':
X\rightarrow \PP_{\RR}^1 \simeq X/\langle\varphi^{g+1}\rangle$ is without real
ramification points, $X$ is given by the real polynomial equation $y^2
=f(x)$, where $f$ is a monic polynomial of degree $2g+2$ and where $f$
has no real roots. By \cite[Prop. 6.3]{G-H} we have $s=1$.
\end{proof}

\begin{rema}
{\rm In the situation of Theorem \ref{maxramnonreel}, we have $\mu
  (\varphi_{\CC}^2 )=4$ and $|\varphi^2|=g+1=g+2-s$. 
  It demonstrates that the inequality
  given in Theorem \ref{boundnonreal} is sharp.}
\end{rema}

\subsection{Automorphisms of prime order}

We bound above the order of an automorphism $\varphi$ of a real curve 
when $|\varphi|$ is prime.
\begin{thm}
\label{primenonreal}
Let $\varphi$ be an automorphism of $X$ of prime order $p$ such that $\pi:
X\rightarrow X/\langle\varphi\rangle$ is without real ramification points. Then
$p\leq g+1$.
\end{thm}

\begin{proof}
Let $g'$ denote the genus of $X'= X/\langle\varphi\rangle$.
Since $p$ is prime, any ramification point $P$ of $\pi_{\CC}$ 
is a totally ramified point 
i.e. $e_P =p$ and $P$ is a fixed point of $\varphi_{\CC}$. 
Consequently, $\pi$ does not have any real branch point i.e. $w_{\RR} (\pi )=0$.

Assume $p\geq g+2$. By Proposition \ref{boundnonrealone} 
we have $\mu (\varphi_{\CC} )\leq 2$.

If $\varphi_{\CC}$ is fixed point free, it follows from (\ref{equ2}) that 
$p=\frac{g-1}{g'-1}$ with $g'\geq 2$. Hence $p\leq g-1$ which gives 
a contradiction.

If $\mu (\varphi_{\CC})=2$ then we get $2g-2=p(2g'-2)+2(p-1)$ by (\ref{equ2}).
Hence $p=\frac{g}{g'}$ and $g'\geq 1$. It gives a contradiction.

We have proved that $p\leq g+1$.
\end{proof}

\subsection{The abelian case}

We give now an upper bound of the order of an abelian group of
automorphisms of a real curve.
\begin{thm}
\label{orderabcurve}
Let $G$ be an abelian group of automorphisms of a real curve $X$ 
such that $\pi :X\rightarrow
X'=X/G$ is without real ramification points.
\begin{description}
\item[({\cal{i}})] If $w (\pi_{\CC} )=0$ then $|G|\leq g-1$.
\item[({\cal{ii}})] If $w (\pi_{\CC} )>0$ then $|G|\leq
  g+3+2(g+1-s)\leq 3g+3$.
\end{description}
\end{thm}

\begin{proof} 
Let $g'$ denote the genus of $X'$. We denote by $s'$ the number
of connected components of $X' (\RR )$.
If $w (\pi_{\CC} )=0$ the Riemann-Hurwitz relation (\ref{equ2}) gives 
$|G|\leq g-1$.

Assume $w_{\CC} (\pi )=0$ and $w_{\RR} (\pi )>0$. Since $\pi :X\rightarrow
X'=X/G$ is without real ramification points, we have $\Int (\pi )=0$
and it follows that 
$s'\geq 2$, hence that $g'\geq 1$ by Harnack inequality. 
By Lemma \ref{pardini}, $w_{\RR} (\pi )>1$. Let $Q_1$, $Q_2$ be two
distinct real branch points of $\pi$ and let $P_1$ be a
point in the fiber $\pi_{\CC}^{-1} (Q_1)$. 
By \cite[Lem. 1.1]{Pa}, $\Stab (P_1)$ is cyclic. Let $\varphi_1$ be a
generator of $\Stab (P_1)$. We have $|\varphi_1|=e(Q_1)$.
Since $G$ is abelian, all the points in the fiber $\pi_{\CC}^{-1}
(Q_1)$ are fixed points of $\varphi_{1,\CC}$. By (\ref{equ6})
and since there are $\frac{|G|}{e(Q_1)}$ points in
the fiber $\pi_{\CC}^{-1} (Q_1)$ , we obtain
\begin{equation}
\label{equ9}
\frac{|G|}{e(Q_1)}\leq 4+2\frac{g+1-s}{e(Q_1) -1}\leq 4+2(g+1-s).
\end{equation}
Similarly, we have $\frac{|G|}{e(Q_2)}\leq 4+2(g+1-s)$.
By (\ref{equ2}) and according to the above remarks, 
we obtain $2g-2\geq (2g' -2)|G|+\sum_{i=1}^2
(|G|-\frac{|G|}{e(Q_i)})\geq 2|G|-8-4(g+1-s)$. Hence $|G|\leq
g+3+2(g+1-s)$.

Assume $w_{\CC} (\pi )>0$. Let $Q_1$, $\bar{Q_1}$ be two conjugate
branch points of $\pi_{\CC}$. Arguing as in the previous case, we get
\begin{equation}
\label{equ8}
\frac{|G|}{e(Q_1)}=\frac{|G|}{e(\bar{Q_1})}\leq 2+(g+1-s).
\end{equation}
If $g'\geq 1$ then, combining (\ref{equ2}) with (\ref{equ8}), we get 
$2g-2\geq (2g' -2)|G|+2
(|G|-\frac{|G|}{e(Q_1)})\geq 2|G|
-4-2(g+1-s)$ i.e. $|G|\leq g+1+(g+1-s)$. 
Let us assume that $g'=0$. Since $s'=1$ and $\pi :X\rightarrow
X'=X/G$ is without real ramification points, we have $w_{\RR} (\pi)=0$.
By (\ref{equ2}), we
see that $w_{\CC} (\pi )\geq 2$ and let $Q_2$, $\bar{Q_2}$ be two conjugate
branch points of $\pi_{\CC}$ distinct from $Q_1$,
$\bar{Q_1}$. Clearly, 
we also have
$\frac{|G|}{e(Q_2)}=\frac{|G|}{e(\bar{Q_2})}\leq 2+(g+1-s)$.
The Riemann-Hurwitz relation (\ref{equ2}) gives
$2g-2\geq  -2|G|+2\sum_{i=1}^2
(|G|-\frac{|G|}{e(Q_i)})\geq 2|G|-8-4(g+1-s)$ i.e. $|G|\leq g+3+2(g+1-s)$. 
\end{proof}

In case of an M-curve, the previous theorem reads:
\begin{cor}
\label{orderabmcurve}
Let $G$ be an abelian group of automorphisms of an M-curve $X$ 
such that $\pi :X\rightarrow
X'=X/G$ is without real ramification points.
\begin{description}
\item[({\cal{i}})] If $w (\pi_{\CC} )=0$ then $|G|\leq g-1$.
\item[({\cal{ii}})] If $w (\pi_{\CC} )>0$ then $|G|\leq
  g+3$.
\end{description}
\end{cor}

\subsection{The hyperelliptic case}

Before giving an upper bound on the order of the automorphisms group
of a real hyperelliptic curve such that the hyperelliptic involution
is without real fixed point, we prove a more general result.
\begin{thm}
\label{ordercenternonrealcurve}
Let $X$ be a real curve such that $\pi :X\rightarrow
X'=X/\Aut (X)$ is without real ramification points.
Assume there exists
$\varphi\in Z(\Aut (X))$, $\varphi\not=Id$, 
such that $\mu (\varphi_{\CC} )>0$, then
$|\Aut (X)|\leq
  2g+2+2\frac{g+1-s}{|\varphi |-1}\leq
  2g+2+2(g+1-s)$.
\end{thm}

\begin{proof} 
Let $g'$ denote the genus of $X'$. We denote by $s'$ the number
of connected components of $X' (\RR )$. Set $N=|\Aut (X)|$.

Let $P$ be a fixed point of $\varphi_{\CC}$, then $\bar{P}$ is also a fixed
point of $\varphi_{\CC}$.

Firstly, assume $Q=\pi_{\CC} (P)$ is a real point. Hence  $Q=\pi_{\CC}
(\bar{P})$.
Since $\pi :X\rightarrow
X'=X/\Aut (X)$ is without real ramification points, we have $\Int (\pi )=0$
and it follows that 
$s'\geq 2$, hence that $g'\geq 1$ by Harnack inequality.
By \cite[Cor. p. 100]{Fa-K},
$\Stab (P)$ is cyclic.
Since the points in the fiber above $Q$
have conjugate stabilizer subgroups and since $\varphi\in Z(\Aut (X))$, 
the fiber
$\pi_{\CC}^{-1} (Q)$ is composed by $\frac{N}{e(Q)}$ points which are 
fixed points of $\varphi_{\CC}$.
By (\ref{equ6}), we obtain
$\frac{N}{e(Q)}\leq \mu (\varphi_{\CC} )\leq 4+2\frac{g+1-s}{|\varphi| -1}\leq
4+2(g+1-s)$.
By (\ref{equ2}) and according to above remarks, 
we obtain $2g-2\geq (2g' -2)N+
(N-\frac{N}{e(Q)})\geq N-\frac{N}{e(Q)}$. 
Hence $N\leq 2g-2+\frac{N}{e(Q)} \leq 2g+2+2\frac{g+1-s}{|\varphi
  |-1}\leq 2g+2 +2(g+1-s)$.

Finally, assume $Q=\pi_{\CC} (P)$ is a non-real point. Hence 
$\bar{Q}=\pi_{\CC} (\bar{P})$ i.e. $P$ and $\bar{P}$ are contained in
two conjugate fibers of $\pi_{\CC}$. Arguing as in the previous case, we get
$
\frac{N}{e(Q)}=\frac{N}{e(\bar{Q})}\leq 2+\frac{g+1-s}{|\varphi
  |-1}$.
If $g'\geq 1$, it follows from (\ref{equ2}) that
$2g-2\geq (2g' -2)N+2
(N-\frac{N}{e(Q})\geq 2N
-4-2\frac{g+1-s}{|\varphi
  |-1}$ i.e. $N\leq g+1+\frac{g+1-s}{|\varphi
  |-1}$. 
Let us assume that $g'=0$.
By (\ref{equ2}), we
see that $w_{\CC} (\pi )\geq 2$ and let $Q_1$, $\bar{Q_1}$ be two conjugate
branch points of $\pi_{\CC}$ distinct from $Q$ and $\bar{Q}$ (the
branch points are
non-real since $g'=0$). Clearly, 
$e(Q_1)=e(\bar{Q_1})\geq 2$.
The Riemann-Hurwitz relation (\ref{equ2}) gives
$2g-2\geq  -2N+2
(N-\frac{N}{e(Q)})+2N(1-\frac{1}{e(Q_1)})\geq N-4-2\frac{g+1-s}{|\varphi
  |-1}$ i.e. 
$N\leq 2g+2+2\frac{g+1-s}{|\varphi
  |-1}$, which completes the proof. 
\end{proof}

Corollary \ref{boundhypreal} and the following proposition  proves that
the order of the automorphisms group of a real hyperelliptic curve
cannot be larger than $4g+4$.

\begin{prop}
\label{orderhypnonrealcurve}
Let $X$ be an hyperelliptic curve such the hyperelliptic involution
$\imath$ does not have any real fixed point.
If $X'=X/\Aut (X)$ is without real ramification points then
then either $g$ is odd and $|\Aut (X)|\leq
4g$, or $|\Aut (X)|\leq 4g+2$.\\
If $X'=X/\Aut (X)$ has at least one real ramification point then
$|\Aut (X)|\leq
4g+4$.
\end{prop}

\begin{proof}
By \cite[Cor. 3 p. 102]{Fa-K}, the hyperelliptic involution $\imath$
is in the center of $\Aut (X)$. Let $g'$ denote the genus of $X'$.
In the case $\pi :X\rightarrow
X'=X/\Aut (X)$ is without real ramification points, the proof follows
from Theorem \ref{ordercenternonrealcurve} and from
\cite[Prop. 4.3]{Mo}.

Now we consider the case when $\pi :X\rightarrow
X'=X/\Aut (X)$ has real ramification points and $\imath$ does not have
real fixed points.
By Proposition \ref{basic1} ({\cal{iv}}), 
there exists two real branch points $Q_1$, $Q_2$ of $\pi$
with real fibers and we have $e (Q_1 )=e(Q_2 )=2$.
Let $P$ be a fixed point of
$\imath$. 
Since the points in the fiber above $Q=\pi_{\CC}(P)$
have conjugate stabilizer subgroups and since $\imath\in Z(\Aut (X))$,
the fiber
$\pi_{\CC}^{-1} (Q)$ is composed by $\frac{N}{e(Q)}$ points which are non-real
fixed points of $\imath$. If $e(Q)=2$ then $\frac{N}{2}\leq \mu
(\imath)=2g+2$
and we get $N\leq 4g+4$. So we assume $e(Q)\geq 4$ ($\Stab (P)$ has
even order since it contains $\imath$).
If $g'\geq 1$, from (\ref{equ2}), we have $2g-2\geq (2g' -2)N+N
(1-\frac{1}{e(Q)})+N
(1-\frac{1}{e(Q_1)})+N(1-\frac{1}{e(Q_2)})\geq
\frac{3N}{4}+\frac{N}{2}+\frac{N}{2}$, hence $N\leq \frac{8}{7}(g-1)\leq
2g-2$. Let us assume that $g'=0$. By (\ref{equ2}), we have $w(\pi_{\CC})>3$.
Let $Q'$ be another branch point of $\pi_{\CC}$ distinct from $Q_1
,Q_2$ and $Q$.
By (\ref{equ2}) and since $\frac{N}{e(Q')}\leq \frac{N}{2}$, we have
$2g-2\geq (2g' -2)N+N
(1-\frac{1}{e(Q)})+N
(1-\frac{1}{e(Q')})+N
(1-\frac{1}{e(Q_1)})+N(1-\frac{1}{e(Q_2)})\geq
-2N+\frac{3N}{4}+\frac{N}{2}+\frac{N}{2}+\frac{N}{2}$, hence $N\leq 8g-8$.
If $ \pi_{\CC}^{-1} (Q)$ is the
unique fiber of $\pi_{\CC}$ containing a fixed point of
$\imath$ then $\frac{N}{e(Q)}= \mu
(\imath)=2g+2$, which contradicts the inequality $N\leq 8g-8$. 
Therefore, we may assume that $\pi_{\CC}^{-1} (Q')$ 
is composed by $\frac{N}{e(Q')}$ points which are non-real
fixed points of $\imath$. Similarly, we assume $e(Q')\geq 4$.
By (\ref{equ2}), we get $2g-2\geq (2g' -2)N+N
(1-\frac{1}{e(Q)})+N
(1-\frac{1}{e(Q')})+N
(1-\frac{1}{e(Q_1)})+N(1-\frac{1}{e(Q_2)})\geq 
-2N+\frac{3N}{4}+\frac{3N}{4}+\frac{N}{2}+\frac{N}{2}$ i.e. 
$N\leq 4g-4$.
\end{proof}


\begin{thebibliography}{ABR2}

\bibitem[Fa-K]{Fa-K} H. M. Farkas, I. Kra,
{\it Riemann Surfaces}, 
Springer-Verlag, 
New York-Berlin-Heidelberg 1980

\bibitem[G-H]{G-H} G. H. Gross, J; Harris, 
{\it Real algebraic curves}, 
Ann. Sci. Ecole Norm. Sup,, 14 (4), 157-182, 1981

\bibitem[Hu]{Hu} J. Huisman, 
{\it On the geometry of algebraic curves having many real components}, 
Rev. Mat. Complut., 14, 83-92, 2001

\bibitem[Kn]{Kn} J. T. Knight, 
{\it Riemman surfaces of field extensions}, 
Proc. Cambridge Philos. Soc., 65, 635-650, 1969

\bibitem[Kr-Ne]{Kr-Ne} W. Krull, J. Neukirch,
{\it Die Struktur der absoluten Galoisgruppe \"{u}ber dem K\"orper
  $\RR (t)$}, 
Math. Ann., 193, 197-209, 1971

\bibitem[Ma1]{Ma1} C. L. May,
{\it Automorphisms of compact Klein surfaces with boundary}, 
Pacific J. Math., 59, 199-210, 1975

\bibitem[Ma2]{Ma2} C. L. May,
{\it Cyclic automorphisms groups of compact bordered Klein surfaces}, 
Houston J. Math., 3, 395-405, 1977

\bibitem[Mo]{Mo} J. P. Monnier, 
{\it Divisors on real curves}, 
Adv. Geom, 3, 339-360, 2003

\bibitem[Pa]{Pa} R. Pardini,
{\it Abelian covers of algebraic varieties}, 
J. Reine angew. Math., 417, 191-213, 1991

\end{thebibliography}
\end{document}